\documentclass[times,sort&compress,3p]{elsarticle}
\journal{Journal of Multivariate Analysis}
\usepackage{amsmath, amssymb, enumerate, amsthm}

\usepackage[usenames,dvipsnames]{xcolor}
\def\tr{{\rm tr\,}} \def\d{\frac{1}{2}\,}
\def\R{\mathbb{R}}\def\E{\mathbb{E}}\def\<{\langle}\def\>{\rangle}
\def\var{{{\mathbb V}\mbox{ar}}}
\def\cov{\mathbb{C}\mathrm{ov}}
\newcommand{\bel}{\begin{equation}\label}
\newcommand{\ee}{\end{equation}}

\title{Multivariate  reciprocal inverse Gaussian  distributions  from the Sabot -Tarr\`es -Zeng  integral}

=-0.5 in
\begin{document}

\begin{frontmatter}

\author[GL]{G\'erard Letac}
\ead{gerard.letac@math.univ-toulouse.fr} 
\author[JW]{Jacek Weso\l owski}
\ead{wesolo@mini.pw.edu.pl}
\address[GL]{Institut de Math\'ematiques de Toulouse, Universit\'e Paul Sabatier, 31062, Toulouse,
France.}
\address[JW]{Matematyki i Nauk Informacyjnych, Politechnika Warszawska, Warszawa, Poland.}  
\begin{abstract} In Sabot and Tarr\`es \cite{ST15}, the authors  have explicitly computed the integral
$$STZ_n=\int \exp( -\<x,y\>)(\det M_x)^{-1/2}dx$$ where  $M_x$ is a symmetric matrix of order $n$ with fixed non-positive off-diagonal coefficients and with  diagonal $(2x_1,\ldots,2x_n)$. The domain of integration is the part of $\R^n$ for which  $M_x$ is positive definite. We calculate more generally for $ b_1\geq 0,\ldots b_n\geq 0$ the integral  $$\int \exp \left(-\<x,y\>-\frac{1}{2}b^*M_x^{-1}b\right)(\det M_x)^{-1/2}dx,$$ we show that it leads to a natural family of distributions in $\R^n$, called the $MRIG_n$ probability laws. This family  is stable by marginalization and by conditioning, and it has number of properties which are multivariate versions
of familiar properties of univariate reciprocal inverse Gaussian
distribution.  In general, if  the power of $\det M_x$ under the integral in $STZ_n$ is distinct from $-1/2$ it is not known how to compute the integral. However, introducing  the graph $G$ having  $V=\{1,\ldots,n\}$ for set of vertices and  the set $E$ of $\{i,j\}'$ s of non-zero entries of $M_x$ as set of edges, we show also that in the particular case where $G$ is a tree,  the integral $$\int \exp( -\<x,y\>)(\det M_x)^{q-1}dx$$ where $q>0,$ is computable in terms of the MacDonald function $K_q.$
\end{abstract}


\


\noindent \textsc{Abbreviated title:} $MRIG_n$ distributions 

\end{frontmatter}

\section{Introduction: the Sabot -Tarr\`es -Zeng integral.}
Let us  describe  the integral appearing  in  Sabot and Tarr\`es \cite{ST15}. Let $W=(w_{ij})_{1\leq i,j\leq n}$
be a symmetric matrix such that $w_{ii}=0$ for all $i=1,\ldots,n$ and such that $w_{ij}\geq 0$ for $i\neq j.$ For $x=(x_1,\ldots,x_n)\in \R^n$ define the matrix $M_x=2\, \mathrm{diag}(x_1,\ldots,x_n)-W.$ For instance if $n=3$ we have
$$M_x=\left[\begin{array}{ccc}2x_1&-w_{12}&-w_{13}\\-w_{12}&2x_2&-w_{23}\\-w_{13}&-w_{23}&2x_3\end{array}\right].$$
Denote by $C_W$ the set of $x\in \R^n$ such that $M_x$ is positive definite. It is easy to see that $C_W$ is an open non-empty unbounded convex set. This is not a cone in general. Frequently we consider  the  undirected graph $G$ with set of vertices $\{1,\ldots,n\}$ and with set of edges
$E=\{\{i,j\}\ ;\ w_{ij}>0 \}$ and we speak of the graph $G$ associated to $W.$
The Sabot-Tarr\`es-Zeng  integral is, for $y_1,\ldots,y_n>0$
\begin{equation}\label{STZI}STZ_n=\int_{C_W}e^{-(x_1y_1+\cdots+x_ny_n)}\frac{dx_1\ldots dx_n}{\sqrt{\det M_x}}=\left(\sqrt{\frac{\pi}{2}}\right)^{n}\frac{1}{\sqrt{y_1\ldots y_n}}e^{-\d\sum_{ij}w_{ij}\sqrt{y_iy_j}}.\end{equation}
 Sabot and Tarr\`es \cite{ST15} give a  probabilistic proof of this remarkable result. Another proof is in Sabot, Tarr\`es and Zeng \cite{STZ17}, based on the Cholesky decomposition. This integral leads naturally to  consideration of probability laws on $\R^n$ that we  call $STZ_n$ distributions with densities proportional to $e^{-\langle x, y\rangle}(\det M_x)^{-1/2}1_{C_W}(x).$ In the present
paper we derive, using a different approach than the two
methods mentioned above,
                a more general  $MRIG_n$ integral in Theorem 2.2. In
particular, we give a new proof of (\ref{STZI}).  The symbol $MRIG$ for multivariate reciprocal inverse Gaussian, is explained below.

This $MRIG_n$  integral enables us to create a new set  (called the $ MRIG_n$ family) of distributions on $\R^n$ which is stable by marginalization and, up to a translation,  stable by conditioning. The bibliography concerning the appearance of the $STZ_n$ and $MRIG_n$ laws in probability theory is already very rich and we suggest to look at Sabot and Zeng \cite{SZ17} and Disertori, Merkl and Rolles \cite{DISERTORIMR17} for many references. An unpublished observation of 2015 of the first author has been used  and reproved in these two publications and some facts of the present paper can be found in them. However, we use here only elementary methods to get our results.

Let us recall that in literature, the generalized inverse Gaussian distributions $ GIG(a,b,q)$ are one dimensional laws with density proportional to $e^{-a^2x-\frac{b^2}{4x}}x^{q-1}1_{(0,\infty)}(x),$  for $a,b>0$ and $q$ real (see Seshadri \cite{SESHADRI} for instance). Parameterizations differ according to the needs of   authors and we have chosen an appropriate one in the present paper. The most  famous particular case is for $q=-1/2$ with the inverse Gaussian distribution.
A  random variable $Y$ with  the  inverse Gaussian distribution $IG(a,b)= GIG(a/2,2b,-1/2)$   has    Laplace transform  \begin{equation}\label{LTIG}\E(e^{-sY})=e^{b(a-\sqrt{a^2+s})}\end{equation}for $s>-a^2$. Its density is proportional to $e^{-\frac{a^2y}{4}-\frac{b^2}{y}}y^{-3/2}1_{(0,\infty)}(y).$   A less known case - but the important one for the present paper- is
for $q=1/2$ with the reciprocal  inverse Gaussian distribution.  Actually, it is a distribution of the inverse of a random variable with an $IG$ distribution.
A  random variable $X$ with  a reciprocal inverse Gaussian distribution $RIG(a,b)= GIG(a,b,1/2)$  has Laplace transform for $s>-a^2$
\begin{equation}\label{LTRIG}\E(e^{-sX})=\frac{a}{\sqrt{a^2+s}}e^{b(a-\sqrt{a^2+s})}\end{equation}  and is such that

\begin{equation}\label{MOMRIG}\E(X)=m=\frac{ab+1}{2a^2}, \ \E(X^2)=\frac{1}{4a^4}(a^2b^2+3ab+3),\ \var(X)=\frac{ab+2}{4a^4}. \end{equation}Its density is proportional to
$e^{-a^2x-\frac{b^2}{4x}}x^{-1/2}1_{(0,\infty)}(x).$  This law is considered for instance in  Barndorff-Nielsen and Koudou \cite{OLEKOUDOU}.
Our  $MRIG_n$
distributions have some properties which are
                multivariate versions of properties known for the
univariate $RIG$ law. These are good reasons for attaching the
name multivariate ($n$-dimensional) $RIG$ to the members of this family. A particular case of the family $MRIG_2$ appears in  Barndorff-Nielsen, Blaesild and Seshadri\cite{OLEBLAESILD}.  The family $STZ_2$ appears in Barndorff-Nielsen and Rysberg \cite{OLERYSBERG}.

 Section 2 proves and comments on the $MRIG_n$ integral, including a presentation of the Disertori-Spencer-Zinbauer\cite{DISERTORISZ10} and Disertori-Merkles-Rolles \cite{DISERTORIMR17} integrals in the studies of supersymmetry.  Section 3 gives some examples. Section 4 details the properties of the $MRIG_n$ laws (we carefully distinguish along the paper the $MRIG_n$ integral and the $MRIG_n$ laws).
 Section 5 considers the particular case of the $STZ_n$ integral when the  graph $G$ associated to $W$  is a tree. Then we  generalize the $STZ_n$ integral by computing in this case
$\int_{C_W} \exp( -\<x,y\>)(\det M_x)^{q-1}dx$ and thus, in particular, obtaining the norming constant for the density considered in Massam and Weso\l owski \cite{MASSAMW04}. Interestingly enough, this generalization allows us to not restrict to the case where the  $w_{ij}$'s are non-negative. The reason is the not so well known fact: if the associated graph of a positive definite matrix $M=(m_{ij})$ is a tree then the symmetric matrix $M'=(\pm m_{ij})$ is still positive definite whatever the $\pm $ are outside of the diagonal; therefore $C_W$ is unchanged.
Section 6  mentions  a striking consequence (Corollary 6.2) of the $MRIG_n$  integral: if $(B_1,\ldots,B_n)$ is multivariate b normal, i.e. $(B_1,\ldots,B_n)\sim N(0,M_x)$, then
$$\Pr(B_1>0,\ldots,B_n>0)=\frac{1}{(2\pi)^{n/2}}\int_{C_W\cap \{t\leq x\}}\frac{dt}{\sqrt{(x_1-t_1)\ldots (x_n-t_n)}\sqrt{\det M_t}}.$$
Section 7 proves a marginal but  delicate fact that  the densities of the $MRIG_n$ distributions are continuous on the whole $\R^n.$ Finally, a first version of this paper is on arXiv 1709.04843.
 \section{The $MRIG_n$ integral}
 \subsection{The integral and its  various forms }
 It is useful to recall  a classical formula, which is in fact  the particular case $n=1$ of  Theorem 2.2 below and the starting point of an induction proof.

\vspace{4mm}\noindent
\textbf{Lemma 2.1.} If $a>0$ and $b\geq 0$ then \begin{equation}\label{BD}\int_0^{\infty}e^{-\frac{a^2t^2}{2}\, -\frac{b^2}{2t^2}}dt=\sqrt{\frac{\pi}{2}}\frac{1}{a}e^{-ab}.\end{equation}

 \vspace{4mm}\noindent

 Various proofs of Lemma 2.1 exist in the literature. An elegant one  considers the equivalent formulation
 \begin{equation}\label{BI}\frac{a}{\sqrt{2\pi}}\int_{-\infty}^{\infty}\exp\left[-\frac{1}{2}(at-\frac{b}{t})^2\right]\, dt =1\end{equation} and proves (\ref{BI}) by the  change of variable $x=\varphi(t)=t-\frac{b}{at}$ which preserves the Lebesgue measure on $\R.$ This idea seems to be due to George Boole\cite{BOOLE}.

 \vspace{4mm}\noindent
\textbf{Theorem 2.2.} Let $a_1,\ldots,a_n>0$ and $b_1,\ldots, b_n\geq 0.$ Then with $a^*=(a_1,\ldots,a_n)$ and $b^*=(b_1,\ldots,b_n)$ we have

\begin{equation}\label{STZI3}MRIG_n=\int_{C_W}e^{-\d\,( a^*M_xa+b^*M_x^{-1}b)}\frac{dx}{\sqrt{\det M_x}}=\left(\frac{\pi}{2}\right)^{n/2}\frac{e^{-(a_1b_1+\cdots+a_nb_n)}}{a_1\ldots a_n}\end{equation}

\vspace{4mm}\noindent
\textbf{Comments.} \begin{itemize}\item Inserting $t=\sqrt{2x}$ in (\ref{BD}) we see that (\ref{BD}) is the particular case $n=1$ of (\ref{STZI3}).
\item Remarkably, the right hand side of (\ref{STZI3}) does not depend on $W.$

\item  Another presentation of (\ref{STZI3}) is
$$\int_{C_W}\exp(-\d\|M_x^{1/2}a-M_x^{-1/2}b\|^2)\frac{dx}{\sqrt{\det M_x}}=\left(\frac{\pi}{2}\right)^{n/2}\frac{1}{a_1\ldots a_n}.$$

For $n=1$ this is nothing but (\ref{BI}) after the change of variable $x=t^2/2.$

\item Another variation: from (\ref{STZI4}) below, writing for short $\sqrt{s}=(\sqrt{s_1},\ldots,\sqrt{s_n})^*$ we have
$$\left(\frac{2}{\pi}\right)^{n/2}\int_{C_W}e^{-\<x,s\>-\frac{1}{2}b^*M_x^{-1}b}\frac{dx}{\sqrt{\det M_x}}=
\frac{1}{\sqrt{s_1\ldots s_n}}e^{-2\<b,\sqrt{s}\>-\frac{1}{2}\sqrt{s}^*W\sqrt{s}}$$

\item One more variation of  (\ref{STZI3}) and  (\ref{STZI4}) is obtained by considering a positive definite matrix $A=(a_{ij})_{1\leq i,j\leq n}$ in the formula

\begin{eqnarray}\nonumber&&\left(\frac{2}{\pi}\right)^{n/2}\int_{C_W}e^{-\d\tr(M_xA)-\frac{1}{2}b^*M_x^{-1}b}\frac{dx}{\sqrt{\det M_x}}\\&=&
\frac{1}{\sqrt{a_{11}\ldots a_{nn}}}e^{-(b_1\sqrt{a_{11}}+\cdots+b_n\sqrt{a_{nn}})-\frac{1}{2}\sum_{i,j=1}^nw_{ij}(\sqrt{a_{ii}a_{jj}}-a_{ij})}\nonumber\end{eqnarray}
If $A=\Sigma^{-1}$, consider the Gaussian random variable $X=(X_1,\ldots,X_n)\sim N(0,\Sigma)$.  Recall that $\rho_{ij}=-a_{ij}/\sqrt{a_{ii}a_{jj}}$ is the correlation between $X_i$ and $X_j$ conditioned by all $(X_k; k\neq i,j).$ Therefore
$$\sum_{i,j=1}^nw_{ij}(\sqrt{a_{ii}a_{jj}}-a_{ij})=\sum_{i,j=1}^nw_{ij}\sqrt{a_{ii}a_{jj}}(1+\rho_{ij}).$$

\item In (\ref{STZI4}) the condition $ a_1,\ldots a_n>0$  is easily relaxed to $ a_1,\ldots a_n\neq 0$: in the right hand side of $ a_1,\ldots a_n>0$  replace $a_i$ by $|a_i|$, Things are quite different for  the condition $ b_1,\ldots b_n\geq 0$: see the comments of the example $n=2$ in Section 3. 

\end{itemize}

\subsection{Proof of Theorem 2.2}

\vspace{4mm}\noindent
\textbf{Proof.} We prove it by induction on $n.$ As mentioned above, Lemma 2.1  is the case $n=1.$ Assume that the result is true for $n.$ Consider
\begin{equation}\label{WM22}W^1=\left[\begin{array}{cc}W&c\\c^*&0\end{array}\right],\ M^1=\left[\begin{array}{cc}M_x&-c\\-c^*&2x_{n+1}\end{array}\right]\end{equation}
where $c=(c_1,\ldots,c_n)^*$ with $c_i\geq 0$ for all $i.$ We now assume that $(x,x_{n+1})\in C_{W^1}.$  From the positive definiteness of $M^1$ we  see that  the Schur complement $t^2=2x_{n+1}-c^*M_x^{-1}c$ is positive. We write 
\begin{equation}\label{MU1}M^1=\left[\begin{array}{cc}I_n&0\\-c^*M_x^{-1}&1\end{array}\right]\left[\begin{array}{cc}M_x&0\\0&t^2\end{array}\right]\left[\begin{array}{cc}I_n&-M_x^{-1}c\\0&1\end{array}\right].\end{equation}  Equality (\ref{MU1}) leads to the computation of $(M^1)^{-1}$ as follows:
\begin{eqnarray}\nonumber(M^1)^{-1}&=&\left[\begin{array}{cc}I_n&M_x^{-1}c\\0&1\end{array}\right]\left[\begin{array}{cc}M_x^{-1}&0\\0&t^{-2}\end{array}\right]\left[\begin{array}{cc}I_n&0\\c^*M_x^{-1}&1\end{array}\right]\\&=&\nonumber
\left[\begin{array}{cc}M_x^{-1}+t^{-2}M_x^{-1}cc^*M_x^{-1}&t^{-2}M_x^{-1}c\\t^{-2}c^*M_x^{-1}&t^{-2}\end{array}\right]\end{eqnarray}
Before writing down the integral $MRIG_{n+1}$ we observe that
\begin{eqnarray*}(a^*,a_{n+1})M^1\left(\begin{array}{c}a\\a_{n+1}\end{array}\right)&=&a^*M_xa-2a^*ca_{n+1}+2x_{n+1}a^2_{n+1}\\&=&-2a^*ca_{n+1}+a^*M_xa+c^*M_x^{-1}ca^2_{n+1}+t^2a^2_{n+1}\end{eqnarray*}
\begin{eqnarray}\nonumber(b^*,b_{n+1})(M^1)^{-1}\left(\begin{array}{c}b\\b_{n+1}\end{array}\right)&=&b^*M_x^{-1}b+t^{-2}b^*M_x^{-1}cc^*M_x^{-1}b+2t^{-2}b^*M_x^{-1}cb_{n+1}+t^{-2}b_{n+1}^2\\&=&b^*M_x^{-1}b+t^{-2}(b_{n+1}+b^*M_x^{-1}c)^2.\label{bbm}\end{eqnarray}

Also observe that the convex set $C_{W^1}$ is parameterized by $(x,t)$ in $C_W\times (0,\infty)$ and that, from (\ref{MU1}) we have $\det M^1=t^2\det M_x .$ With this parameterization we have $$\frac{dxdx_{n+1}}{\sqrt{\det M^1}}=\frac{dx}{\sqrt{\det M_x}}dt.$$We now write $MRIG_{n+1}$ as follows

\begin{eqnarray}\nonumber  MRIG_{n+1}&=&\nonumber e^{a^*ca_{n+1}}\int_{C_W}\exp-\frac{1}{2}\left[a^*M_xa+c^*M_x^{-1}ca^2_{n+1}+b^*M_x^{-1}b\right]\\&&\nonumber\left(\int_0^{\infty}\exp-\frac{1}{2}\left[t^2a^2_{n+1}+t^{-2}(b_{n+1}+b^*M_x^{-1}c)^2\right]dt\right)\frac{dx}{\sqrt{\det M_x}}\\&=&\nonumber \sqrt{\frac{\pi}{2}}\frac{1}{a_{n+1}}e^{a^*ca_{n+1}-a_{n+1}b_{n+1}}\\&&\label{MU4}\int_{C_W}\exp-\frac{1}{2}\left[a^*M_xa+(c^*a_{n+1}+b^*)M_x^{-1}(ca_{n+1}+b)\right]\frac{dx}{\sqrt{\det M_x}}\\&=&\label{MU5}\left(\frac{\pi}{2}\right)^{(n+1)/2}\frac{1}{a_1\ldots a_{n+1}}e^{-a^*b-a_{n+1}b_{n+1}}
\end{eqnarray}
In this chain of equalities (\ref{MU4}) is a consequence of  Lemma 2.1 applied to the pair $$a_{n+1} ,b_{n+1}+b^*M_x^{-1}c.$$
Here a comment is in order: a famous lemma of Stieltjes implies that $M_x^{-1}$ has non-negative coefficients when $x\in C_W.$ Let us detail  the  proof in this particular case: if $D=2\, \mathrm{diag}(x_1,\ldots,x_n)$ then $M_x=D^{1/2}(I_n-A)D^{1/2}$ where $A=D^{-1/2}WD^{-1/2}$. Since $M_x$ is positive definite, $I_n-A$ is also positive definite. Now write $(I_n-A)^{-1}=I_n+A+\ldots+A^{2N-1}+A^N(I_n-A)^{-1}A^N.$
Since $A^N(I_n-A)^{-1}A^N$  is  positive semidefinite, its trace is $\geq 0$ and therefore for all $N$
$$\sum_{k=0}^{2N-1}\tr(A^k)\leq \tr(I_n-A)^{-1}$$
Since $A$ has non-negative coefficients this implies that $\sum_{k=0}^{\infty}\tr(A^k)$ converges.
In particular $\lim_{N\rightarrow \infty}\tr(A^{2N})=0.$ This implies that all the eigenvalues of $A$ are in $(-1,1)$ and therefore the series of matrices $S=\sum_{k=0}^{\infty}A^k$ converges to $(I_n-A)^{-1}.$ Since $A$ has non-negative coefficients the same is true for $S$ and for $M_x^{-1}=D^{-1/2}SD^{-1/2}.$ Furthermore, if the graph $G$ has  vertices $\{1,\ldots,n\}$ and has edges $\{i,j\}$  present according to the fact that $a_{ij}>0$ or not, then $(I_n-A)^{-1}$ is positive definite if $G$ is connected (this remark will be used in the proof of Lemma 2.4 below).

As a consequence $b_{n+1}+b^*M_x^{-1}c\geq 0$ and therefore (\ref{BD}) is applicable.
 Equality (\ref{MU5}) is a consequence of the induction hypothesis where the pair $(a,b)$ is replaced by $(a,  a_{n+1}\, c+b).$ The induction hypothesis is extended. $\square$

\subsection{ Laplacian of $W$  and  parameterizations of $C_W$ by $(0,\infty)^n$ and $\R^n$} In order to show in Section 2.4 that two other remarkable integrals can be deduced from the  $MRIG_n$ integral \eqref{STZI3}, it is necessary to recall some definitions about Laplacian on graphs or weighted graphs (see for instance Bapat ]cite{BAPAT}).

We define the Laplacian of $W$ as the quadratic form on $\R^n$ defined by
\begin{equation}\label{LAPLACIAN}v^*L_Wv=\sum_{i<j}w_{ij}(v_i-v_j)^2=\d \sum_{i=1}^n\sum_{j=1}^nw_{ij}(v_i-v_j)^2.\end{equation} If $s_i=\sum_{j=1}^nw_{ij}$ and if $D=\mathrm{diag}(s_1,\ldots,s_n)$ the representative matrix of this quadratic form is $L_W=D-W.$ From the definition it is semi positive definite, and since $(1,\ldots,1)^*$ is a eigenvector of $L_W$ associated to the the eigenvalue zero, $L_W$ cannot be positive definite.  However, by adding a proper diagonal matrix $$D_b=\mathrm{diag}(b_1,\ldots,b_n)$$ with $b_i\geq 0$ the matrix $D_b+L_W$ can be positive definite. One can also  remark that $W^1=\left[\begin{array}{cc}W&b\\b^*&0\end{array}\right]$ implies that $ L_{W^1}=\left[\begin{array}{cc}D_b+L_W&-b\\-b^*&\sum_{j=1}^nb_j\end{array}\right].$

\vspace{4mm}\noindent\textbf{Lemma 2.3.} $D_b+L_W$ is positive definite if and only if
 for each connected component $C$ of the graph  associated to $W$ there exists $k\in C$  such that $b_k>0.$

\vspace{4mm}\noindent\textbf{Proof.} $\Leftarrow$ Enough is to assume that  the associated graph is connected and that there exists a $k$ such that $b_k>0.$ If $v$ is such that $v^*(D_b+L_W)v=0$ then $v_k=0.$ Furthermore $v_i-v_j=0$ if $w_{ij}>0$. Since the associated graph is connected all the $v_i$'s are equal, and they are zero like $v_k:$ this shows the positive definiteness of $D_b+L_W.$ $\Rightarrow$ Here again we can assume that $G$ is connected. We have seen that if $b_i=0$ for all $i$ then $D_b+L_W=L_W$ cannot be positive definite. $\square$

 The next lemma describes an important parameterization of $C_W$ by $(0,\infty)^n.$ Note that it depends on a non-zero parameter $b\in [0,\infty)^n.$ The case $b=(1,\ldots,1)^*$ is most useful.

\vspace{4mm}\noindent\textbf{Lemma 2.4.} Assume that the graph $G$ associated to $W$  is connected.  Let $y\in (0,\infty)^n$, fix $b\in [0,\infty)^n$ such that $b\neq 0$ and define $x\in \R^n$ by
\begin{equation}\label{XVT}2x_i=\frac{1}{y_i}\left(b_i+\sum_{j=1}^nw_{ij}y_j\right)\end{equation}
Then $x$ belongs to $C_W$, we have $M_x=D_bD^{-1}_y+L_W$ and $y=M^{-1}_x b,$ the map $y\mapsto x$ is a diffeomorphism from  $(0,\infty)^n$ onto $C_W$
and  \begin{equation}\label{XVY}dx=\frac{\det{M_x}}{2^n}\frac{\, dy}{y_1\ldots y_n}.\end{equation}

\vspace{4mm}\noindent\textbf{Proof.} We rewrite \eqref{XVT} as $2x_iy_i-\sum_{j=1}^nw_{ij}y_j=b_i$ and thus it is equivalent to $b=M_xy.$ Denote \begin{equation}\label{WY}W^{(y)}=D_yWD_y\end{equation} for a while and
$s_i^{(y)}=\sum_{j=1}^nw_{ij}y_iy_j=2x_iy_i^2-b_iy_i.$ Therefore

$$D_{s^{(y)}}=2D_yD_xD_y-D_bD_y,\ \ \ L_{W^{(y)}}=D_{s^{(y)}}-W^{(y)}
=2D_yD_xD_y-D_bD_y-D_yWD_y.$$ From the definition \eqref{LAPLACIAN} of the Laplacian we have $ L_{W^{(y)}}=D_yL_WD_y$ and we get
$$D_bD_y +L_{W^{(y)}}=2D_yD_xD_y-D_yWD_y,\ \  D_bD^{-1}_y+L_W=M_x.$$
From Lemma 2.3 $M_x=D_bD^{-1}_y+L_W$ is positive definite and furthermore  $y=M^{-1}_x b.$ Equality $b=M_xy$ shows that the map $y\mapsto x $ from $(0,\infty)^n$ to $C_W$ is injective since $0=M_x(y-y')$ implies $y=y'$ from the definite positiveness of $M_x.$  If $x\in C_W$, define    $y=M^{-1}_x b.$ The fact that $M_x^{-1}$ has only non-negative coefficients implies that   $y\in [0,\infty)^n.$  The fact that $G$ is connected implies that  $y\in (0,\infty)^n.$ We get   $M_x=D_bD^{-1}_y+L_W$ and this shows the surjectivity since any $y\in (0,\infty)$ provides a positive definite matrix $D_bD^{-1}_y+L_W$.  The fact that $y\mapsto x$  is a diffeomorphism from $(0,\infty)^n$ onto $C_W$ is clear.

The differential of the map $y\mapsto x$  from $C_W$ onto  $(0,\infty)^n$ is
\begin{equation}\label{DIFFY}h\mapsto -2M_x^{-1}D_hM_x^{-1}b=-2M_x^{-1}D_hy\end{equation}
For showing \eqref{DIFFY} we observe that the differential of $M\mapsto M^{-1}$ is $H\mapsto -M^{-1}HM^{-1}$  and that the differential of the map $x\mapsto M_x$ is $h\mapsto 2\, D_h.$ The Jacobian of $y\mapsto x$ is therefore $\frac{2^n}{\det M_x}y_1\cdots y_n$ and this proves \eqref{XVY}. $\square$

Replacing $y_i$ by $e^{t_i}$ we will use Lemma 2.4 in the next section under the following form:

\vspace{4mm}\noindent\textbf{Corollary 2.5.} Under the hypothesis of Lemma 2.4, for $t\in \R^n$ define
$$2x_i(t)=b_ie^{-t_i}+\sum_{j=1}^nw_{ij}e^{t_j-t_i}.$$ Then the map $t\mapsto x=x(t)$ is a diffeomorphism from $\R^n$ onto $C_W$ and $dx=\frac{\det M_x}{2^n}dt.$

\subsection{The Disertori-Spencer-Zirnbauer and Disertori-Merkl-Rolles integrals } In application of Theorem 2.2 and Corollary 2.5, we prove two surprizing formulas $DSZ_n$ and $DMR_n$ due to Disertori, Spencer and Zirnbauer \cite{DISERTORISZ00} and  Disertori, Merkl and Rolles  \cite{DISERTORIMR15} (2015).  For describing them we need the following notations. We consider the quadratic form in (1.1) of the first paper:

$$v^*D(t)v=\sum _{1\leq i<j\leq n}w_{ij}e^{t_i+t_j}(v_i-v_j)^2+\sum_{k=1}^nb_ke^{t_k}v_k^2.$$ The element $(i,i)$ of the corresponding $n\times n$ matrix $D(t)$ is
$b_ie^{t_i}+\sum_{j=1}^nw_{ij}e^{t_i+t_j}$ and the off diagonal element $(i,j)$ is $-w_{ij} e^{t_i+t_j}.$ This is nothing but the quadratic form with matrix
$D(t)=D_bD_y+L_{W^{(y)}}$ as in  \eqref{WY} when $y_i=e^{t_i}$ for all $i. $

We introduce a function $G(t)$ which is only marginally   different from the $F$ defined by (1.2) in\cite{DISERTORISZ00}.
\begin{equation}\label{FG}G(t)=\sum_{i<j}w_{ij}(\cosh(t_i-t_j)-1)+\sum_{k=1}^n\left((\cosh t_k-1)b_k+t_k\right).\end{equation}With these notations,  the surprising formula (1.4) of \cite{DISERTORISZ00},  see \eqref{DSZI4} below, is the subject of the following proposition.

\vspace{4mm}\noindent\textbf{Proposition 2.6.} Assume that $W$ is such that the associated graph is connected  and fix  $b\in [0,\infty)^n$ with $b\neq 0.$ Then
\begin{equation}\label{DSZI4}DSZ_n=\frac{1}{(\sqrt{2\pi})^n}\int_{\R^n}e^{-G(t)} \sqrt{\det D(t)}dt=1.\end{equation}

\vspace{4mm}\noindent \textbf{Proof.} In \eqref{STZI3} we insert $a_1=\ldots=a_n=1$  and we make the change of variable $x\mapsto t$ from $C_W$ onto $\R^n$ described in Corollary 2.5. We get \begin{eqnarray*}-\d(a^*M_{x(t)}a+b^*M_{x(t)}^{-1}b)&=&-\d \sum_{i=1}^n2x_i(t)+\d\sum_{i=1}^n\sum_{j=1}^nw_{ij}-\d\sum_{i=1}^nb_ie^{t_i}\\&=&-\d \sum_{i=1}^nb_i(e^{t_i}+e^{-t_i})-\d\sum _{i<j}w_{ij}(e^{t_j-t_i}+e^{t_i-t_j}-2)\\&=& -\sum_{i=1}^nb_i\cosh t_i-
\sum _{i<j}w_{ij}(\cosh (t_j-t_i)-1).
\end{eqnarray*}
Since $D(t)=D_{y(t)}M_{x(t)}D_{y(t)}$ we have $\det D(t)=e^{\sum_{i=1}^n 2t_i}\det M_{x(t)}.$  Finally using Corollary 2.5 we obtain \eqref{DSZI4}. $\square$

\vspace{4mm}\noindent Similarly, formula (2.4) of Disertori, Merkl and Rolles\cite{DISERTORIMR15} 
introduces a  probability $\mu(ds,dt)$ on $\R^n\times \R^n$ defined by
\begin{equation}\label{MU}\mu(ds,dt)=e^{-\d s^*D(t)s-G_1(t)}\det{D(t)}\frac{dtds}{(2\pi)^n}\end{equation} where the function $G_1$ is quite close to the function  $G$ defined by \eqref{FG} and is defined by
$$G_1(t)=\sum_{i<j}w_{ij}(\cosh(t_i-t_j)-1)+\sum_{k=1}^n (e^{-t_k}b_k+t_k)=G(t)+\sum_{k=1}^n(1-\sinh t_k)b_k.$$ If $(S,T)\sim \mu$ it is clear that $S$ is Gaussian when conditioned by $T$. However, the fact that  the total mass of  $\mu$ is  one is not  that obvious. The result is stated in Proposition 2.8 below. We skip its proof which uses again Corollary 2.5. It is  a consequence of the $STZ_n$ integral \eqref{STZI}, by  the change of variable of Corollary 2.5. The hypotheses on $W$ and $b$ are the same  as in Proposition 2.6.

\vspace{4mm}\noindent\textbf {Proposition 2.7} Let $f(t)=\frac{1}{(\sqrt{2\pi})^n}e^{-G_1(t)} \sqrt{\det D(t)}.$ Then $f$ is a probability density on $\R^n$. Furthermore if $T\sim f(t)dt$ and $S|T\sim N(0,D(T)^{-1})$ then $(S,T)\sim \mu$ defined by \eqref{MU}.

\section{Examples} The following examples consider various graphs associated to $W$ where some  calculations about $M_x$ are  explicit.
\subsection {The case $n=2.$} We take $$W=\left[\begin{array}{cc}0&1\\1&0\end{array}\right],\ M_x= \left[\begin{array}{cc}2x_1&-1\\-1&2x_2\end{array}\right],\ M_x^{-1}= \frac{1}{4x_1x_2-1}\left[\begin{array}{cc}2x_2&1\\1&2x_1\end{array}\right]$$ and $C_W$ is the convex set of $\R^2$ limited by one branch of a hyperbola
$$C_W=\{(x_1,x_2)\ ; \ x_1,x_2>0,\ 4x_1x_2-1>0\}.$$ Theorem 2.2 says that
$$\int\int_{C_W}\exp-[a_1^2x_1+a_2^2x_2-a_1a_2+\frac{1}{4x_1x_2-1}(b_1^2x_2+b_2^2x_1+b_1b_2)]\frac{dx_1dx_2}{\sqrt{4x_1x_2-1}}=\frac{\pi}{2}\frac{e^{-a_1b_1-a_2b_2}}{a_1a_2}.$$
A warning: the extension of $MRIG_n$ to the case where some $b_i$' s are negative leads to a non elementary elementary integral. The case $n=2$ is appropriate for explaining this fact: following the steps of the proof of Theorem 2.2 we arrive  up to a multiplicative constant to the integral 
$$e^{a_1a_2}\int_0^{\infty}e^{-a_1^2x_1-\frac{a_2^2+b_1^2}{4x_1}-a_2\left|b_2+\frac{b_1}{2x_1}\right|}\frac{dx_1}{\sqrt{2x_1}}$$ that we cannot evaluate when $b_1b_2<0.$

\subsection {The complete graph for $n\geq 3$}
We consider the case where $w_{ij}=c$ for all $i\neq j.$ Denote by $J_n$ the $n\times n$ matrix with all entries equal to $1.$ Therefore $W=c(J_n-I_n).$

\vspace{4mm}\noindent
\textbf{Proposition 3.1.} If $W=c(J_n-I_n)$  then \begin{equation}\label{EM1}\det M_x=(c+2x_1)\ldots (c+2x_n)(1-\sum_{i=1}^n\frac{c}{c+2x_i})\end{equation}
\begin{equation}\label{EM2}C_W=\{(x_1,\ldots,x_n); x_1,\ldots,x_n\geq 0 ,\ \sum_{i=1}^n\frac{c}{c+2x_i}<1\}\end{equation}

\vspace{4mm}\noindent
\textbf{Proof.} Write $D=2\mathrm{diag}(x_1,\ldots,x_n)+cI_n$. Therefore $M_x=D-cJ_n=D^{1/2}(I_n-A)D^{1/2}$ where $A=cD^{-1/2}J_nD^{-1/2}=vv^*$ and $$v=(\frac{\sqrt{c}}{\sqrt{c+2x_1}},\ldots,\frac{\sqrt{c}}{\sqrt{c+2x_n}})^*.$$

The eigenvalues of $A$ are 0 with multiplicity $(n-1)$ and $v^*v=\sum_{i=1}^n\frac{c}{c+2x_i}.$ This implies that the eigenvalues of $I_n-A$ are $1$ with multiplicity $n-1$ and
$1-\sum_{i=1}^n\frac{c}{c+2x_i}$. This leads to (\ref{EM1}). To prove  (\ref{EM2}), clearly the right-hand side contains  $C_W$. Conversely if $x_i>0$ for all $i$ writing $M_x=D^{1/2}(I-A)D^{1/2}$ shows $x\in C_W$ if and only if $I-A$ is positive definite, ie
$1-\sum_{i=1}^n\frac{c}{c+2x_i}>0.$  $\ \ \square$
\subsection {The daisy}
 We consider the case where
 $$W=\left[\begin{array}{ccccc}0&c_1&c_2&\ldots&c_n\\c_1&0&0&\ldots&0\\c_2&0&0&\ldots&0\\ \ldots&\ldots&\ldots&\ldots&\ldots\\c_n&0&0&\ldots&0\end{array}\right],\ M_x=\left[\begin{array}{ccccc}2x_0&-c_1&-c_2&\ldots&-c_n\\-c_1&2x_1&0&\ldots&0\\-c_2&0&2x_2&\ldots&0\\ \ldots&\ldots&\ldots&\ldots&\ldots\\-c_n&0&0&\ldots&2x_n\end{array}\right]$$
It is easy to see by induction that
$$\det M_x=2^nx_1\ldots x_n\left(2x_0-\sum_{i=1}^n\frac{c_i^2}{2x_i}\right),$$$$ C_W=\{(x_0,\ldots,x_n); x_0,\ldots,x_n>0,\  2x_0-\sum_{i=1}^n\frac{c_i^2}{2x_i}>0\}.$$ It is elementary to write $M_x^{-1}$ explicitly. If we write for simplicity
$$M(a,b,c)=\left[\begin{array}{ccccc}a&c_1&c_2&\ldots&c_n\\c_1&b_1&0&\ldots&0\\c_2&0&b_2&\ldots&0\\ \ldots&\ldots&\ldots&\ldots&\ldots\\c_n&0&0&\ldots&b_n\end{array}\right]$$  $B=b_1\ldots b_n$, $D=\det M(a,b,c)=B\left(a-\sum_{i=1}^n\frac{c_i^2}{b_i}\right)$ and $H=M(a,b,c)^{-1}=(h_{ij})_{0\leq i,j\leq n}$ then for $  i\neq j$ and distinct from $0$ we have
$$h_{00}=\frac{B}{D},\  h_{0i}=-\frac{B}{D}\frac{c_i}{b_i},\ h_{ii}=\frac{1}{b_i}\left(1+\frac{B}{D}\frac{c_i^2}{b_i}
\right),\  h_{ij}=\frac{c_ic_j}{b_ib_j}\left(1+\frac{B}{D}(\frac{c_i^2}{b_i}+\frac{c_j^2}{b_j})
\right)$$ For $n=2$ it gives $D=\det M_x=8x_0x_1x_2-2x_2c_1^2-2x_1c_2^2$ and

$$M_x^{-1}=\frac{1}{D}\left[\begin{array}{ccc}4x_1x_2&2c_1x_2&2c_2x_1\\2c_1x_2&4x_0x_2-c_2^2&c_1c_2\\2c_2x_1&c_1c_2&4x_0x_1-c_1^2\end{array}\right]$$
\subsection{The chain} Define  $A_{n+1}$ as the graph $\stackrel{0}{\bullet}-\stackrel{1}{\bullet}-\stackrel{2}{\bullet}-\ldots-\stackrel{n}{\bullet}$ corresponding to the matrix

$$W=\left[\begin{array}{ccccccc}0&c_1&0&0&\cdots&0&0
\\ c_1&0&c_2&0&\cdots&0&0
\\  0 &c_2&0&c_3&\cdots&0&0
\\ 0&0&c_3&0&\cdots&0&0
\\ \cdots &\cdots&\cdots&\cdots&\cdots&\cdots&\cdots
\\ 0&0&0&0&\cdots&0&c_n
\\    0 &0&0&0&\cdots&c_n&0 \end{array}\right]$$
where $c_1,\ldots,c_n>0.$ Thus $M_x=-W+\mathrm{diag}(2x_0,2x_1,\ldots,2x_n)$ is a Jacobi matrix. Without losing  generality we may assume that $c_1=\ldots=c_n=1$ by the transformation $$\mathrm{diag}(\lambda_0,\lambda_1,\ldots,\lambda_n)\ M_x\ \mathrm{diag}(\lambda_0,\lambda_1,\ldots,\lambda_n)$$ where
$$\lambda_0=1,\  \lambda_1=\frac{1}{c_1 },\ \lambda_{2p}=\frac{c_1 c_3\ldots c_{2p-1}}{c_2c_4\ldots c_{2p}},\ \lambda_{2p+1}=\frac{c_2c_4\ldots c_{2p}}{c_1 c_3\ldots c_{2p+1}}$$ thus replacing $x_i$ by the affinity $x_i\lambda^2_i.$ If $D_0=2x_0$ and $  D_1=4x_0x_1-1$ then the determinant $D_n$ of the matrix
$$M_x=\left[\begin{array}{ccccccc}2x_0&-1&0&0&\cdots&0&0
\\ -1&2x_1&-1&0&\cdots&0&0
\\  0 &-1&2x_2&-1&\cdots&0&0
\\ 0&0&-1&2x_3&\cdots&0&0
\\ \cdots &\cdots&\cdots&\cdots&\cdots&\cdots&\cdots
\\ 0&0&0&0&\cdots&2x_{n-1}&-1
\\    0 &0&0&0&\cdots&-1&2x_n \end{array}\right]$$
is computable by the induction formula $D_n=2x_nD_{n-1}-D_{n-2}.$ For instance for $n=3$ the set  $C_W$ is described by the four inequalities
$$x_0>0,\  4x_0x_1-1>0,\ 4x_0x_1x_2-x_0-x_2>0, \  16x_0x_1x_2x_3-4x_0x_3-4x_2x_3-4x_1x_2+1>0.$$
Since the chain $A_{n+1}$ is also a tree, results of Section 5 below are applicable to this example.

\section{A study of the $MRIG_n$ distributions}
\subsection{$MRIG_n$  as a natural exponential family}   Writing $s_i=a_i^2$ in  \eqref{STZI3}  we obtain 

\begin{equation}\left(\frac{2}{\pi}\right)^{n/2}\int_{C_W}e^{-\<x,s\>-\frac{1}{2}b^*M_x^{-1}b}\frac{dx}{\sqrt{\det M_x}}=
\frac{1}{\sqrt{s_1\ldots s_n}}e^{-(b_1\sqrt{s_1}+\cdots+b_n\sqrt{s_n})-\frac{1}{2}\sum_{i,j=1}^nw_{ij}\sqrt{s_is_j}}\label{STZI4}\end{equation}
which suggests that the natural exponential family (NEF)  concentrated on $C_W\subset \R^n$  generated by the unbounded measure
\begin{equation}\label{EXPbW}\mu(b,W)(dx)=e^{-\frac{1}{2}b^*M_x^{-1}b}1_{C_W}(x)\frac{dx}{\sqrt{\det M_x}}\end{equation}
is interesting to study. For $n=1$  if $b_1>0$ this is nothing but a $RIG$ distribution mentioned in (\ref{LTRIG})  and if $b_1=0$ it is  a Gamma family with shape parameter $1/2$. For $n>1$ and $b=0$ this family is considered in Sabot, Tarr\`es and Zeng (2016).
Given $W$ and $a\in (0,+\infty)^n,\ b\in [0,+\infty)^n$ we consider the probability on $[0,+\infty)^n$ defined by
$$P(a;b,W)(dx)=(\frac{2}{\pi})^{n/2}\left(\prod_{j=1}^na_je^{a_jb_j}\right)e^{-\d a^*M_xa-\d b^*M_x^{-1}b}\textbf{1}_{C_W}(x)\frac{dx_1\ldots dx_n}{\sqrt{\det M_x}}.$$ We say that $P(a;b,W)$ is a $ MRIG_n$ distribution. Theorem 2.2 proves that it is indeed a probability. From time to time we will use the notation $f(a;b,W)(x)$ for the density of $P(a;b,W).$  Note that  $(X_1,\ldots,X_n)\sim P(a;b,0)$ iff $X_1,\ldots,X_n$ are independent and $X_k\sim RIG(a_k,\,b_k)$, $k=1,\ldots,n$.

 In this section,   we show that if $X$  has a $MRIG_n$ distribution then the subvector $(X_1,\ldots,X_k)$ has a $MRIG_k$ distribution. We also show that up to a translation factor, the conditional distribution of $(X_{k+1},\ldots,X_n)$ given $(X_1,\ldots,X_k)$  has  a $MRIG_{n-k}$ distribution. Thus the class of  $MRIG_n$  distributions has a remarkable property of stability by marginalization and conditioning. These facts have been independently observed by Sabot and Zeng (2019) in their Lemma 5, and also mentioned in Sabot and Zeng \cite{SZ17} quoting the arXiv versions of  Sabot and Zeng \cite{SZ17} 
and of the present paper.

We begin with the calculation of  the Laplace transform of $P(a;b,W).$ Introducing the following function:
\begin{eqnarray} \label{FLT}G(a;b,W)&=&(\frac{\pi}{2})^{n/2}\left(\prod_{j=1}^n\frac{e^{-a_jb_j}}{a_j}\right) e^{-\d a^*Wa}\end{eqnarray}
we remark that $P(a;b,W)(dx)$ can be written as
\begin{equation}\label{EXPb}P(a;b,W)(dx)=\frac{1}{G(a;b,W)}e^{-x_1a_1^2-\ldots-x_na_n^2}\mu(b,W)(dx)\end{equation}
Under the form (\ref{EXPb}) we see that for fixed $b\in [0,\infty)^n$
$$F_b=\{P(a;b,W), a\in (0,\infty)^n \}$$ is a natural exponential family, parameterized by $a$ and not by its natural parameter $(s_1,\ldots,s_n)=(a^2_1,\ldots,a_n^2) ,$  and generated by  $\mu(b,W)$.  From the fact that the mass of (\ref{EXPb}) is one, the Laplace transform of $\mu(b,W)$ is
defined for $s\in (0,\infty)^n$ by$$L_{\mu(b,W)}(s)=G((\sqrt{s_1},\ldots,\sqrt{s_n});b,W).$$ We deduce from this the  form of the Laplace transform of $P(a;b,W)$ itself.

\vspace{4mm}\noindent \textbf{Proposition 4.1.} If $(X_1,\ldots,X_n)\sim P(a;b,W)$ then
\begin{equation}\label{LAPLACE}\E(e^{-s_1X_1-\ldots-s_nX_n})=e^{\<a,b\>-\<\sqrt{a^2+s},b\>}e^{a^*Wa-\sqrt{a^2+s}^*W\sqrt{a^2+s}}\prod_{j=1}^n\frac{a_j}{\sqrt{a_j^2+s_j}}\end{equation}
where we have written symbolically $\sqrt{a^2+s}=(\sqrt{a_1^2+s_1}, \ldots,\sqrt{a_n^2+s_n})^*.$ In particular
\begin{eqnarray}\label{ESP}\E(X_i)=m_i&=&\frac{1}{2a_i}\left(b_i+\sum_{j\neq i}w_{ij}a_j\right), \\ \label{VAR}\var(X_i)&=&\frac{1}{4a_i^4}+\frac{m_i}{2a_i^2}\\\label{COV}
\cov(X_i,X_j)&=&-\frac{w_{ij}}{4a_ia_j}\end{eqnarray}

\vspace{4mm}\noindent \textbf{Comments.}  The one dimensional margins are the classical $RIG$ distributions (\ref{LTRIG}). More specifically the distribution of $X_i$ is $$RIG(a_i,b_i+\sum_{j=1}^nw_{ij}a_j)=RIG(a_i,2a_im_i).$$ In other terms the two  parameters  of the distribution of $X_i$'s are the $i$ components of the vectors $a$ and $b+Wa.$

Formula \eqref{ESP} expresses $m_i$  with a formula which is the successful change of variable \eqref{XVT} where the pair $(x,y)$ is replaced here by $(m,a).$ Observe also that the covariance of $(X_i,X_j)$  is never positive.

It can be mentioned that, like for a Gaussian distribution, the parameters $(a,b,W)$ of the distribution $MRIG_n$ are determined if we know the distributions of all pairs $(X_i,X_j): $ the knowledge of the distribution of $X_i$ gives from  \eqref{ESP} and \eqref{VAR} the knowledge of $ m_i$ and $a_i.$ The knowledge of the distribution of $(X_i,X_j)$ and of $a$ gives from \eqref{COV} the knowledge of $w_{ij}$ and $W,$ and then \eqref{ESP} gives the value of $b_i.$ Estimation of the parameters can be designed from this remark.
One more analogy with the Gaussian distributions is the fact that if $X\sim MRIG_n$ then $X_i$ and $X_j$ are independent if and only if they are uncorrelated: this can be read from the Laplace transform of $X.$

\vspace{4mm}\noindent \textbf{Proof of Proposition 4.1.} Formula \eqref{LAPLACE} comes immediately from $$\E(e^{-s_1X_1-\ldots-s_nX_n})=\frac{G(\sqrt{a_1^2+s_1}, \ldots,\sqrt{a_n^2+s_n}); b,W)}{G(a;b,W)}.$$ Formulas \eqref{ESP} and \eqref{VAR} are consequence of the properties of the one dimensional $RIG$ given in \eqref{MOMRIG}. The simple formula \eqref{COV} is obtained by $$\cov(X_i,X_j)=\frac{\partial^2}{\partial s_i \partial s_j}\log G(\sqrt{a_1^2+s_1}, \ldots,\sqrt{a_n^2+s_n}); b,W)\vert_{s=0}.\ \square$$

\subsection{The marginals of the $ MRIG_n$ distribution}

For stating the next results we need the following notations:
\begin{itemize}

\item For   vectors $(a_1,a_2,\ldots,a_n)^*$ and $(b_1,b_2,\ldots,b_n)^*$ we denote
$$\tilde{a}_k=(a_1,\ldots,a_k)^*,\ \tilde{b}_k=(b_1,\ldots,b_k)^*.$$ With this notation sometimes we write $P(\tilde{a}_n;\tilde{b}_n,W)$ for $P(a;b,W).$

\item If $$W_n=\left[\begin{array}{ccccc}0&w_{12}&w_{13}&\ldots&w_{1n}\\
w_{12}&0&w_{23}&\ldots&w_{2n}\\
\ldots&\ldots&\dots&\ldots&\ldots\\
w_{1n}&w_{2n}&w_{3n}&\ldots&0\end{array}\right]$$ for $k=2,3,\ldots,n$ we take $c_k\in \R^{k-1}$ defined as the $k$th column of $W_n$ but restricted to be above the diagonal, namely $c_k=(w_{1k},w_{2k},\ldots,w_{k-1,k})^*.$
\item If $k=1,2,3,\ldots,n$ we write $W_n$ by blocks as follows
$$W_n=\left[\begin{array}{cc}W_k&W_{k,n-k}\\ \ W_{k,n-k}^*&W'_{n-k}\end{array}\right]$$
In other terms $W_k=[w_{ij}]_{1\leq i,j\leq k},\  W'_{n-k}=[w_{ij}]_{k+1\leq i,j\leq n}.$
\item We define the killing symbol $K$ from $\R^n$ to $\R^{n-1}$ by $$K(x_1,\ldots,x_n)^*=(x_1,\ldots,x_{n-1})^*$$ For instance $K \tilde{b}_k= \tilde{b}_{k-1}.$ In general  for $k<n$ we have $$K^{n-k}(x_1,\ldots,x_n)^*=(x_1,\ldots,x_k)^*$$

\end{itemize}

\vspace{4mm}\noindent \textbf{Proposition 4.2.} If $(X_1,\ldots,X_n)\sim P(\tilde{a}_n,\tilde{b}_n,W)$ then $
(X_1,\ldots,X_k)\sim P(\tilde{a}_k,B_k,W_k)$ where
$$B_k=\tilde{b}_k+\sum_{j=k+1}^na_jK^{j-k-1}c_j=\tilde{b}_k+W_{k,n-k}(a_{k+1},\ldots,a_n)^*$$

\vspace{4mm}\noindent \textbf{Proof.} For $k=n-1$, this is claiming that $
(X_1,\ldots,X_{n-1})\sim P(\tilde{a}_{n-1},\tilde{b}_{n-1}+a_nc_n,W_{n-1}).$ Such a formula is essentially formula (\ref{MU4}) when replacing $n$ by $n+1.$

We now prove the result by induction on $n-k.$ Suppose that  $
(X_1,\ldots,X_k)\sim P(\tilde{a}_k,B_k,W_k)$  is true. Then as for the passage from $n$ to $n-1$ we can claim that  $
(X_1,\ldots,X_{k-1})\sim P(\tilde{a}_{k-1},KB_k+a_kc_k,W_{k-1}).$  Now we have
\begin{eqnarray*}KB_k+a_kc_k&=&K\tilde{b}_k+a_kc_k+K\sum_{j=k+1}^na_jK^{j-k-1}c_j\\&=&\tilde{b}_{k-1}+a_kc_k+\sum_{j=k+1}^na_jK^{j-k}c_j=B_{k-1}\end{eqnarray*} and the induction is extended. $\square$

\vspace{4mm}\noindent \textbf{Comments.}\begin{itemize}
\item Proposition 4.2 could have been proved with the Laplace transform of Proposition 4.1, but is seems that after all induction is simpler.
\item A reformulation of Proposition 4.2 is the explicit form of the integral $$\int_{C_{W'_{n-k}}}f_{\tilde{a}_n;\tilde{b}_n,W_n}(\tilde{x}_k,x_{k+1},\ldots,x_n)dx_{k+1}\ldots dx_n=f_{\tilde{a}_k, B_k,W_k}(\tilde{x}_k),$$
namely

\begin{eqnarray}\nonumber&&(\frac{2}{\pi})^{n/2}\left(\prod_{j=1}^na_j\right)e^{\langle a, b\rangle-\d a^*Wa}\int_{C_{W'_{n-k}}}
e^{-(x_1a^2_1+\cdots+x_na^2_n)-\d b^*M_x^{-1}b}\textbf{1}_{C_W}(x)\frac{dx_{k+1}\ldots dx_n}{\sqrt{\det M_x}}\\&& =\nonumber(\frac{2}{\pi})^{k/2}\left(\prod_{j=1}^ka_j\right)e^{\langle \tilde{a}_k,B_k\rangle}e^{-\d \tilde{a}_k^*M_{\tilde{x}_k}\tilde{a}_k-\d B_k^*M_{\tilde{x}_k}^{-1}B_k}\textbf{1}_{C_{W_k}}(\tilde{x}_k)\frac{1}{\sqrt{\det M_{\tilde{x}_k}}}.
\end{eqnarray}

\item  Inserting $b=0$ in Proposition 4.2 makes that $(X_1,\ldots,X_n)$ has an $STZ_n$ distribution. If we also take $k=n-1$ we see that $B_{n-1}=a_nc$ where $c=(w_{i,n})_{i=1}^{n-1}.$ As a consequence, we see that any $MRIG_{n-1}$ distribution is a projection of some $STZ_n$ distribution.  This explains why Sabot, Tarr\`es and Zeng \cite{STZ15} indeed observe that  one dimensional  margins of an $STZ_n$ distribution are $RIG$ ones.

\end{itemize}

\subsection{Conditional distributions under $MRIG_n$} Let us begin by some general observations about exponential families on a product  $E\times F$ of two Euclidean spaces generated by the distribution $\pi(dx)K(x,dy)$. Let $\Theta\subset E\times F$ be the interior of the set
$$\{(t,s) \ ; \; L(t,s)=\int_{E\times F}e^{-\<t,x\>-\<s,y\>}\pi(dx)K(x,dy)<\infty\}.$$ We assume that $\Theta$ is the product of two open subsets of $E$ and $F$ respectively:
\begin{equation}\label{HYPI}\Theta=\Theta_E\times \Theta_F\end{equation}
We fix $(t_0,s_0)\in \Theta$ and we consider a random variable $(X,Y)$ valued in $E\times F$ with density
$$\frac{1}{L(t_0,s_0)}e^{-\<t_0,x\>-\<s_0,y\>}\pi(dx)K(x,dy).$$ We are interested in the Laplace transform of the conditional distribution of $Y|X.$ For computing this, we consider the marginal density of $X$ with respect to $\pi:$
$$\frac{1}{L(t_0,s_0)}\int_Fe^{-\<t_0,x\>-\<s_0,y\>}K(x,dy)= e^{-\<t_0,x\>}\frac{g(s_0;x)}{L(t_0,s_0)}$$
where we have introduced the auxiliary function
$$g(s_0;x)=\int_F e^{-\<s_0,y\>}K(x,dy)$$ defined on $\Theta_F\times E.$ As a consequence, the conditional distribution of $Y|X$ is $e^{-\<s_0,y\>}K(X,dy)/g(s_0;X)$ and its Laplace transform is for $s+s_0\in \Theta_F$ the ratio
\begin{equation}\label{LTCY}s\mapsto \frac{g(s+s_0;X)}{g(s_0;X)}\end{equation} Suppose now that we are able to identify a density on $F$ having (\ref{LTCY}) as Laplace transform. In this case the problem of computation of the density of $Y|X$ will be solved.

We are going to apply this program to $E=\R^k$, $F=\R^{n-k}$, to a probability  on $\R^n$ defined by its density  $$f(x)\propto e^{-\d b^*M_x^{-1}b}\textbf{1}_{C_W}(x)\frac1{\sqrt{\det M_x}}$$ and finally to $t_0=(a_1^2,\ldots,a_k^2),\ s_0=(a_{k+1}^2,\ldots,a_n^2).$ We prove in Section 7 that $f$ is continuous on $\R^n$ when $b_1,\ldots,b_n$ is not zero and when the graph associated to $W$ is connected.

We have that condition (\ref{HYPI}) is fulfilled with $\Theta_E=(0,\infty)^k$ and $\Theta_F=(0,\infty)^{n-k}.$ Of course $\tilde{X}_k=(X_1,\ldots,X_k)$ and $(X_{k+1},\ldots,X_n)$ replace $X$ and $Y$. Also $s$ is now $(s_{k+1},\ldots.s_n)$ and $s+s_0$ is described by
$$A_k(s)=(\sqrt{a_{k+1}^2+s_{k+1}},\ldots, \sqrt{a_{n}^2+s_{n}})^*.$$ We now deduce the crucial function $g(s_0;x)$ from Proposition 4.2 where the marginal law of $(X_1,\ldots,X_k)$  is computed. We get
\begin{equation}\label{TGF}g(a_{k+1}^2,\ldots,a_n^2; \tilde{x}_k)=\frac{G(a;b,W)}{G(\tilde{a}_k; B_k,W_k)}\frac{e^{-\d B^*_kM^{-1}_{\tilde{x}_k}B_k}}{\sqrt{\det(M_{\tilde{x}_k})}}\textbf{1}_{C_{W_k}}(\tilde{x}_k)\end{equation} where the function $G(a;b,W)$ has been introduced in (\ref{FLT}).
Remember that the right hand side of (\ref{TGF})  depends  of $a_{k+1},\ldots,a_n$ also through $B_k=\tilde{b}_k+W_{k,n-k}(a_{k+1},\ldots,a_n)^*.$ Let us adopt the notation
\begin{equation}\label{Bs}
B_k(s)=\tilde{b}_k+W_{k,n-k}A_k(s)\end{equation}
Here is now the Laplace transform of $(X_{k+1},\ldots,X_n)$  given $\tilde{X}_k$
\begin{eqnarray}\nonumber\E(e^{-s_{k+1}X_{k+1}-\ldots-s_nX_n}|\tilde{X}_k=\tilde{x}_k)&=&\frac{g((a_{k+1}^2+s_{k+1},\ldots,a_n^2+s_n); \tilde{x}_k)}{g(a_{k+1}^2,\ldots,a_n^2; \tilde{x}_k)}=PQR\\\label{IN1}P &=& \frac{G((\tilde{a}_k,A_k(s));b,W)}{G(a;b,W)}\\\label{IN2}Q&=&\frac{G(\tilde{a}_k;B_k,W_k)}{G(\tilde{a}_k;B_k(s),W_k)}\\\label{IN3}R&=&
e^{-\d B_k(s)^*M^{-1}_{\tilde{x}_k}B_k(s)+\d B_k^*M^{-1}_{\tilde{x}_k}B_k}
\end{eqnarray}
It is our intention to prove the existence of $\alpha=(\alpha_{k+1}, \ldots,\alpha_n)^*,$ $\beta=(\beta_{k+1}, \ldots,\beta_n)^*,$  $\gamma=(\gamma_{k+1}, \ldots,\gamma_n)^*$
 and of a matrix $\mathcal{W}$ such that $$PQR=\frac{G(\sqrt{\alpha_{k+1}^2+s_{k+1}},\ldots,\sqrt{\alpha_{n}^2+s_{n}}; \beta, \mathcal{W})} 
{G({\alpha_{k+1},\ldots,\alpha_n}; \beta, \mathcal{W})}
e^{-\gamma_{k+1}s_{k+1}-\ldots -\gamma_ns_n}$$ which is
saying that the conditional distribution of $(X_{k+1}-\gamma_{k+1},\ldots,X_n-\gamma_n)$ given $\tilde{X}_k$ is a $MRIG_{n-k}$ distribution.
The next proposition gives the complete result:

\vspace{4mm}\noindent \textbf{Proposition 4.3.} For $X\sim P(a;b,W)$ consider $\alpha =( a_{k+1},\ldots,a_n)^*$, $\beta\in \R^{n-k}$ defined by
$$\beta=(b_{k+1},\ldots,b_n)^*+W_{n-k,k}M^{-1}_{\tilde{X}_k}\tilde{b}_k,$$ the matrix $D=\mathrm{diag}(\gamma_{k+1},\ldots,\gamma_n)$ defined as the diagonal part of the matrix $W_{n-k,k}M^{-1}_{\tilde{X}_k}W_{k,n-k}$ and
$$\mathcal{W}=W'_{n-k}+W_{n-k,k}M^{-1}_{\tilde{X}_k}W_{k,n-k}-D.$$
Then  the conditional distribution  of $(X_{k+1}-\gamma_{k+1},\ldots,X_n-\gamma_n)$  given $\tilde{X}_k$ is    $ P(\alpha,\beta,\mathcal{W}).$

\vspace{4mm}\noindent \textbf{Proof.}  We have to analyze the dependency on $s$ of the three quantities $P,Q,R$ defined above by (\ref{IN1}), (\ref{IN2}), (\ref{IN3}).
However for simplification, we do not write the factors which do not depend on $s$. More specifically we introduce the following   equivalence relation among non-zero functions $f$ or $g$ depending on $s$ and possibly on other parameters like $a$, $b$, $W$ by writing $f\equiv g$ if $f(s)/g(s)$ does not depend on $s.$ For instance
\begin{eqnarray*}P&\equiv &\prod_{j=k+1}^n(a^2_j+s_j)^{-\d} \times e^{-b_{k+1}\sqrt{a_{k+1}^2+s_{k+1}}-\cdots-b_{n}\sqrt{a_{n}^2+s_n}}e^{-\tilde{a}_k^*W_{k,n-k}A_k(s)-\d A_k(s)^*W'_{n-k}A_k(s)}\\
Q&\equiv& e^{\tilde{a}_k^*W_{k,n-k}A_k(s)}\\
R&\equiv &e^{-\d B_k(s)^*M^{-1}_{\tilde{X}_k}B_k(s)}\equiv e^{-A_k(s)^*W_{n-k,k}M^{-1}_{\tilde{X}_k}\tilde{b}_k}\times e^{-\d A_k(s)^*W_{n-k,k}M^{-1}_{\tilde{X}_k}W_{k,n-k}A_k(s)}
\end{eqnarray*}
 A patient analysis of the product $PQR$ as a function of $A_k(s)$ gives Proposition 4.3.
$\square$

\subsection{A convolution property of the $MRIG_n$ laws}

The following  proposition is a generalization of the following additive convolution:
$$RIG(a,b)*IG(a,b')=RIG(a,b+b'),$$ (see Barndorff-Nielsen and Koudou \cite{OLEKOUDOU}, Barndorff-Nielsen and Rydberg \cite{OLERYDBERG} and Barndorff-Nielsen, Blaesild and Seshadri \cite{OLEBLAESILD})
with definitions in (\ref{LTIG}) and (\ref{LTRIG}).

\vspace{4mm}\noindent \textbf{Proposition 4.4.} Let  $a_i,b_i, b'_i>0$ for $i=1,\ldots,n.$ If $X=(X_1,\ldots,X_n)$ has the $MRIG_n$ distribution $ P(a;b,W)$, if $Y=(Y_1,\ldots,Y_n)$   such that  $Y_i\sim IG(a_i,b_i')$ with independent components,  and if $X$ and $Y$ are independent, then
$$X+Y=(X_1+Y_1,\ldots,X_n+Y_n)\sim P(a; b+b',W).$$

\vspace{4mm}\noindent \textbf{Proof.} Just compute the Laplace transform, using Proposition 4.1 and (\ref{LTIG}). $\square$

\subsection{Questions}Here are some unsolved problems linked to $MRIG_n$ laws

 \begin{itemize}\item If $X\sim P(a;b,W)$ what is the distribution of $M^{-1}_X?$ This random matrix is concentrated on a manifold of dimension $n.$ This is a natural question since in one dimension if $X\sim RIG(a,b)$ then the distribution of $1/X$ is known and is $IG(b/2,2a).$ However the Laplace transform of $M^{-1}_X$, namely $L(s)=\mathbb{E}(e^{-\tr(sM^{-1}_X)})$ defined when $s$ is a positive definite matrix of order $n,$ is not known in general. If $b=0$ then  $X$ has an $STZ_n$ distribution and Theorem 2.2 shows that $L(s)$ is known  for $s$ of rank one.

\item Since in one dimension $IG$ and $RIG$ distributions are particular cases of the generalized inverse Gaussian laws, the natural extension of the $MRIG_n$ laws  is to consider
the probability densities on $\R^n$ proportional to
\begin{equation}\label {MRIGq}e^{-\d a^*M_xa-\d b^*M_x^{-1}b}(\det M_x)^{q-1}\bf{1}_{C_W}(x)\end{equation}
extending our familiar $MRIG_n$ integral from $1/2$ to an arbitrary real number $q$. But the corresponding integral extending Theorem 2.2 is untractable. However, in a particular case, namely  if $b=0$ and if the graph $G$ associated to $W$ is a tree, Proposition  5.1  below computes the integral on $C_W$ of the function (\ref{MRIGq}).  A related distribution has been analyzed in Massam and Weso\l owski \cite{MASSAMW04}in connection with a multivariate version of the Matsumoto-Yor property (see  e.g. Matsumoto and Yor\cite{MYOR}, Letac and Weso\l owski \cite{LETACW} and Massam and Weso\l owski \cite{MASSAMW06}).

\item  Probabilistic interpretations of the one dimensional laws $IG$ and $RIG$ are known, as hitting time and time of last visit of an interval $[a,\infty)$ by a drifted
Brownian motion $ t\mapsto mt+B(t)$ (in the respective cases $m>0$ and $m<0$). How to extend this to $MRIG_n$ laws? An answer to this problem is given in Sabot and Zeng (2017) but one may look for other interpretations.

\end{itemize}

\section{Another generalization of the Sabot-Tarr\`es-Zeng  integral: the case of a tree}

In this section we consider another generalization of a specialization of the Sabot-Tarr\`es-Zeng  integral (\ref{STZI}): we assume that the graph $G$ associated to $W$  is a tree but we replace in (\ref{STZI}) the power $-1/2$ of $\det M_x$ by the real number $q-1>-1$. Furthermore in Proposition 5.2, we are able to drop the restriction $w_{ij}\geq 0$ that we have done all along the paper, because of the following proposition of linear algebra:

\vspace{4mm} \noindent\textbf{Proposition 5.1.} Let $M=(m_{ij})_{1\leq i,j\leq n}$ be a symmetric matrix and let
$$E=\{(i,j); 1\leq i<j\leq n, \ m_{ij}\neq 0.\}$$
Assume that $G$ is a graph with set of vertices $\{1,\ldots,n\}$ and with $E$ as set of edges. Then \begin{enumerate}
\item If $G$ is a tree or a forest, $\det M$ is a polynomial in $(m_{ii})_{i=1}^n$ and in $(m_{ij}^2)_{(i,j)\in E}.$ \item If $G$ is a tree or a forest,  if  $M$ is positive definite and if $(\epsilon_{ij})_{1\leq i,j\leq n}$ is a symmetric matrix such that $\epsilon_{ij}=\pm 1$ and $\epsilon_{ii} =1$ for all $i,j $ then the symmetric matrix $(\epsilon_{ij}m_{ij})_{1\leq i,j\leq n}$ is also positive definite. \item If the graph has a cycle then $\det M$ is a sum of monomials such that at least one of them contains an odd power of some $m_{ij}$ with ${i,j}\in E$.  \end{enumerate}

\vspace{4mm} \noindent\textbf{Comments.} In general, changing the two off-diagonal entries $(ij)$ and $(ji)$ of a positive definite matrix $M$  into their  opposite creates a new symmetric matrix which is not positive definite anymore. The proposition shows that this is not the case when the graph associated to $M$ is a tree. Part 3 shows that the fact that $\det M$
 is a polynomial in $(m_{ij}^2)_{I,J}\in E$ characterizes the fact that the graph is a tree or a forest.

\vspace{4mm} \noindent\textbf{Proof.} We prove 1)  by induction on $n.$ The result is clear when $n=1$ and $n=2.$ Suppose that it is true for $n$ and consider the case of a symmetric matrix $M_1$ of order $n+1$  such that its associated graph $G_1$ is a tree. Without loss of generality, we may assume that $n+1$ has only one neighbour in the tree and that this neighbor is $n.$ This implies that $M_1$ has the form
$$M_1=\left[\begin{array}{cc}M&v\\v^*&m_{n+1,n+1}\end{array}\right],\ \ v^*=(0,\ldots,0,m_{n+1,n})$$ where the symmetric matrix $M$ of order $n$  is associated to  the graph $G$ which is $G_1$ minus the vertex $n+1.$ Since $n+1$ had $n$ as the only neighbour, $G$ is also a tree. Write also
$$M=\left[\begin{array}{cc}M_{-1}&v_{-1}\\v_{-1}^*&m_{n,n}\end{array}\right]$$
where $M_{-1}$ is symmetric of order $n-1$. Assume that $m_{n+1,n+1}\neq 0$ and denote $c=m_{n+1,n}^2/m_{n+1,n+1}$ and
$$M'=\left[\begin{array}{cc}M_{-1}&v_{-1}\\v_{-1}^*&m_{n,n}-c\end{array}\right]$$
Therefore we get that
$$\det M_1=m_{n+1,n+1}\det M'.$$ Since $M'$ is symmetric and since its associated graph is the tree $G$, the induction hypothesis implies that $\det M'$ is a polynomial with respect to the squares of the $m_{ij}$ with $1\leq i<j\leq n$ where $(i,j)$ is an edge of $G$. Also,  $\det M'$ is an affine function of $c$. Since $c$  is a multiple of $m_{n+1,n}^2$, therefore the extension of the induction hypothesis to $n+1$ is done when  $m_{n+1,n+1}\neq 0.$ The extension to the case $m_{n+1,n+1}=0 $ is done by continuity of the polynomial $\det M_1.$

For showing 2) we now apply 1) to the case where $M$ is positive definite and we assume without loss of generality that $G$ is a tree.  We number its vertices $\{1,\ldots,n\}$ such that if $G_k$ is the graph associated to the restriction $M_k$ of $M$ to $\{1,\ldots,k\}^2$, then $G_k$ is a tree, a point which can be proved by induction. Denote
$ M_k(\epsilon)=(\epsilon_{ij}m_{ij})_{1\leq i,j\leq k}.$ Since $G_k$ is a tree and since $M$ is positive definite, then $\det(M_k(\epsilon))>0$. From the theorem of principal determinants, $M_n(\epsilon)$ is positive definite.

For showing 3) we assume first that $G$ contains the cycle $1-2-\ldots-n-1.$ We choose $m_{ij}=0$ if $|i-j|\neq 1$ $m_{12}=m_{21}=a$ and $m_{i,j}=1$ for the other edges of the cycle. With this choice the matrix $M$ is $$M_n=\left[\begin{array}{ccccccc}0&a&0&0&\cdots&0&1\\ a&0&1&0&\cdots&0&0\\  0 &1&0&1&\cdots&0&0\\ 0&0&1&0&\cdots&0&0\\ \cdots &\cdots&\cdots&\cdots&\cdots&\cdots&\cdots\\ 0&0&0&0&\cdots&0&1\\    1 &0&0&0&\cdots&1&0 \end{array}\right].$$
Standard techniques show that $\det M_n=\det M_{n+4}$ for $n\geq 3$ and that
$$\det M_3=\det M_5=2a, \ \det M_4=(a-1)^2, \ \det M_6=-(a+1)^2.$$
Therefore one of the monomials is $\pm 2a$: this odd power is the one which  was announced and this ends the proof of Proposition 5.1. $\square$

\vspace{4mm} \noindent
For stating Proposition 5.2 we need to introduce the MacDonald function on $ (0,\infty):$
$$K_q(x)=\d\int_0^{\infty }u^{q-1}e^{-\d x(u+\frac{1}{u})}du$$ It is useful to display a property of this integral  \begin{equation}
\label{MCDO2}
2\left(\frac{b}{a}\right)^qK_q(2ab)=\int_0^{\infty}v^{q-1}e^{-a^2v-\frac{b^2}{v}}dv
\end{equation}
We denote by  $s(i)$  the number of neighbours of $i$ in the tree, namely the size of $\{j: w_{ij}>0\}.$

\vspace{4mm} \noindent\textbf{Proposition 5.2.} Let  $W=(w_{ij})_{1\leq i,j\leq j}$ be  a symmetric matrix with zero diagonal such that its associated graph $G$ is a tree. Let $M_x=2\, \mathrm{diag}(x_1,\ldots,x_n)-W$ and let $C_W$ be the set of $x$'s such that $ M_x$ is positive definite. If $q>0$ then
 \begin{equation}\label{MWF}\int_{C_W}e^{-\d a^*M_xa}(\det M_x)^{q-1}dx=
2^{q-1}\Gamma(q)e^{\d a^*Wa}\prod_{i=1}^na_i^{q(s(i)-2)}\prod_{i<j}|w_{ij}|^qK_q(a_ia_j|w_{ij}|)\end{equation}

\vspace{4mm} \noindent\textbf{Comments.}\begin{itemize}\item For $y_1,\ldots,y_n> 0$ another presentation  of \eqref{MWF} is
$$\int_{C_ W}e^{-\langle x,y \rangle}(\det M_x)^{q-1}dx=
2^{q-1}\Gamma(q)\prod_{i=1}^ny_i^{\d q(s(i)-2)}\prod_{i<j}w_{ij}^qK_q(\sqrt{y_iy_j}w_{ij}).$$
\item Of course, inserting  $q=1/2$ gives back the Sabot -Tarr\`es -Zeng- integral in the case where $G$ is a tree. To check this we use Lemma 2.1 above which says
$$
K_{1/2}(x)=\sqrt{\tfrac{\pi}{2x}}\,e^{-x},\quad x>0.
$$
 For $q=3/2$ we use  Watson \cite{WATSON} page 90 formula 12 for getting  $$
K_{3/2}(x)=\sqrt{\tfrac{\pi}{2x}}\,e^{-x}\left(1+\frac{1}{x}\right), \quad x>0.
$$  and we obtain
$$\int_{C_W}e^{-\d a^*M_xa}\sqrt{\det M_x}dx=\left(\frac{\pi}{2}\right)^{n/2}\prod_{i=1}^na_i^{-3}\prod_{i<j}(1+a_ia_jw_{ij}).$$
\item We give a proof of Proposition 5.1, while another proof could be extracted from Massam and Weso\l owski (2004), where the authors consider the NEF generated by the unbounded measure
$$1_{C_W}(x)(\det M_x)^{q-1}dx$$ and  independence properties of  distributions from this NEF.
 Bobecka \cite{BOBECKA} has a multivariate generalization.

\end{itemize}

\vspace{4mm} \noindent\textbf{Proof.} We proceed by induction on $n$. This is correct for $n=1$ since in this case  $s(1)=0$ and since the empty product $\prod_{i<j}$ is one. Suppose that the formula (\ref{MWF}) is true for $n$ and let us extend  it to $n+1.$  We use the same notation as in Section 2.2:  we keep the notations $a,$  $W$ and  $M_x$ for the matrices of order $n$ as before and we consider the block matrices $M^1$ and $W^1$ defined by  (\ref{WM22}). We now use a different factorization of $M^1$ by writing
\begin{equation}\label{CW1}M^1=
\left[\begin{array}{cc}I_n&-\frac{c}{2x_{n+1}}\\0&1\end{array}\right]
\left[\begin{array}{cc}M_x-\frac{cc^*}{2x_{n+1}}&0\\0&2x_{n+1}\end{array}\right]\left[\begin{array}{cc}I_n&0\\-\frac{c^*}{2x_{n+1}}&1\end{array}\right].\end{equation}
Since the graph $G^1$ which is associated to $W^1$ is  a tree, without loss of generality we assume that the vertex $n+1$ has only one neighbour which is $n.$
In other terms, we may assume that the vector $c$ of $\R^n$ has the form
$$c=(0,\ldots,0,w_{n,n+1})^*.$$ This choice implies also that the graph $G$ associated to $W$ is still a tree. Formula (\ref{CW1}) implies that $C_{W^1}$ is the set of $(x,x_{n+1})\in \R^{n+1}$ such that $x_{n+1}>0$ and such that in the matrix
$$M_y=M_x-\tfrac{cc^*}{2x_{n+1}}=M_x-
\left[\begin{array}{cc}0&0\\0&\frac{w_{n,n+1}^2}{2x_{n+1}}\end{array}\right]$$ its half diagonal
$y=(x_1,\ldots,x_{n-1},x_n-\frac{w_{n,n+1}^2}{4x_{n+1}})^*$
belongs to $C_W.$ The Jacobian of the transformation $(x,x_{n+1})\mapsto (y,x_{n+1})$ is one. Therefore we can write
\begin{eqnarray*}&&\int_{C_{W^1}}
e^{-\d a^*M_xa+a^*ca_{n+1}-a_{n+1}^2x_{n+1}}(\det M^1)^{q-1}dxdx_{n+1}\\
&=&e^{ a^*ca_{n+1}}\int_{C_{W^1}}
e^{-\d a^*(M_x-\frac{cc^*}{2x_{n+1}})a-\frac{(c^*a)^2}{4x_{n+1}}-a_{n+1}^2x_{n+1}}\det (M_x-\frac{cc^*}{2x_{n+1}})^{q-1}(2x_{n+1})^{q-1}dxdx_{n+1}\\
&=&\int_{C_{W}}
e^{-\d a^*M_y a}
(\det M_y)^{q-1}dy\times e^{ a^*ca_{n+1}}\int_0^{\infty}e^{-\frac{(c^*a)^2}{4x_{n+1}}-a_{n+1}^2x_{n+1}}(2x_{n+1})^{q-1}dx_{n+1}\\
\end{eqnarray*}
The latter integral is expressed with (\ref{MCDO2}) as
$$e^{|w_{n,n+1}|a_na_{n+1}}\frac{a_n^q}{a_{n+1}^q}|w_{n,n+1}|^qK_q(a_na_{n+1}|w_{n,n+1})|,$$ and the former one is (\ref{MWF}), from the induction hypothesis. To conclude, observe that the number of neighbours of $n+1$ in $G^1$ is one, and that the number of neighbours of $n$ in $G^1$ is the number $ s(n)$ of neighbours of $n$ in $G$ plus one. $\square$

\section{If $B\sim N(0,M_x)$ what is $\Pr(B_1>0,\ldots,B_n>0)$?} Of course the exact solution of this question cannot be found. However  for any $x\in C_W$ denote by $$f(x)=\Pr(B_1>0,\ldots,B_n>0).$$ Then formula
(\ref{STZI4}) enables us to compute the Laplace transform of $f(x)\mathbf{1}_{C_W}(x).$
The trick is to observe that the first member of (\ref{STZI4}) involves the density $g_x(b)$ of $N(0,M_x)$ when the $b_1,\ldots,b_n$ are restricted to be $>0.$ Of course  $f(x)=\int_{\R^n_+}g_x(b)db$ , and $b\mapsto g_x(b)/f(x)$ is a probability density on $\R^n_+$.

One can even get a knowledge of the Laplace transform of $b\mapsto g_x(b).$ More specifically

\vspace{4mm} \noindent\textbf{Proposition 6.1.} For $\theta=(\theta_1,\ldots,\theta_n)\in \R^n_+$ and $x\in C_W,$ denote
$\int_{\R^n_+}e^{-\<\theta,b\>}g_x(b)db=f(x,\theta).$ Then for $y=(y_1,\ldots,y_n)\in \R^n_+$ we have
\begin{equation}\label{HHH}\int_{C_W}e^{-\<x,y\>}f(x,\theta)dx=\frac{1}{2^n\sqrt{y_1(y_1+\theta_1)}\ldots\sqrt{y_n(y_n+\theta_n)}}e^{-\frac{1}{2}\sum_{i,j=1}^nw_{ij}\sqrt{y_iy_j}}.\end{equation}In particular

$$\int_{C_W}e^{-\<x,y\>}f(x)dx=\frac{1}{2^ny_1\ldots y_n}e^{-\frac{1}{2}\sum_{i,j=1}^nw_{ij}\sqrt{y_iy_j}}.$$

\vspace{4mm} \noindent\textbf{Proof.} Enough is to multiply both sides of (\ref{STZI4})
 by $e^{-\<\theta,b\>}$ and integrate with respect to  $b$ on $\R^n_+.$ Permuting the integrations on the left hand side leads to (\ref{HHH}). $\square$

\vspace{4mm} \noindent\textbf{Corollary 6.2.}
\begin{equation}\label{FOR}f(x)=\frac{1}{(2\pi)^{n/2}}\int_{C_W\cap \{t\leq x\}}\frac{dt}{\sqrt{(x_1-t_1)\ldots (x_n-t_n)}\sqrt{\det M_t}}\end{equation}

\vspace{4mm} \noindent\textbf{Proof.} Denote $$h(x)=\frac{1}{\pi^{n/2}(x_1\ldots x_n)^{1/2}}\mathbf{1}_{\R^n_+}(x), \ \ g(x)=\frac{1}{(2\pi)^{n/2}}\frac{1}{\sqrt{\det M_x}}\mathbf{1}_{C_W}(x).$$  Consider the Laplace transforms $L_f(y)$, $L_g(y),$ $L_h(y)$ defined for $y_1,\ldots,y_n>0$. They are given respectively  by (\ref{HHH}) with $\theta=0$ , by the Sabot-Tarr\`es-Zeng  integral  (\ref{STZI}) and by $$\int_{\R^n_+}e^{-\<x,y\>}h(x)=\frac{1}{\sqrt{y_1\ldots y_n}}.$$ As a consequence $L_f=L_gL_h$ which implies that $f$ is the convolution product of $g$ and $h$ and proves (\ref{FOR}). $\square$

\vspace{4mm} \noindent\textbf{Corollary 6.3.} With the notation $D_y=\mathrm{diag}(y_1,\ldots,y_n)$ and $y=M^{-1}_xb$ where $b=(1,\ldots,1)^*$ we have
\begin{equation}\label{MONST}f(x)=\frac{1}{(2\pi)^{n/2}}\sqrt{y_1\ldots y_n}\int_{(0,\infty)^n}\left(\frac{\det(D_{u+y}^{-1}+L_W)}{\prod_{i=1}^n(u_i(u_i+y_i)))}\right)^{1/2}du_1\ldots du_n
\end{equation}
with the Laplacian $L_W$ defined in \eqref{LAPLACIAN}.

\vspace{4mm} \noindent\textbf{Proof.} In \eqref{FOR} we make the change of variable introduced in Lemma 2.4, namely $s=M^{-1}_tb.$ We also observe that  again from Lemma 2.4 it follows that $M_x-M_t=D_b(D_y^{-1}-D_s^{-1})$ and that
$$\prod_{i=1}^n(x_i-t_i)=\frac{1}{2^n}\det(M_x-M_t)=\frac{1}{2^n}\det D_b(D_y^{-1}-D_s^{-1})=\frac{1}{2^n}\prod_{i=1}^nb_i\frac{s_i-y_i}{s_iy_i}.$$
Using the fact that $b=(1,\ldots,1)^*$ we get
$$f(x)=\frac{1}{(2\pi)^{n/2}}\int_{y_1}^{\infty}\ldots\int_{y_n}^{\infty}\left(\det (D_s^{-1}+L_W)\prod_{i=1}^n\frac{y_i}{s_i(s_i-y_i)}\right)^{1/2}ds_1\ldots ds_n.$$  Hence the change of variables: $u_i=s_i-y_i$, $i=1,\ldots,n$, yields  \eqref{MONST}.$\square$

\vspace{4mm} \noindent\textbf{Comments.}
\begin{itemize}\item Applying formula (\ref{FOR}) even to the case $n=2$ is surprising, since the left hand side is explicitly known: recall that if $$(X_1,X_2) \sim N\left(0,\left[\begin{array}{cc}1&-\cos \alpha\\-\cos \alpha&1\end{array}\right]\right) \\ \Rightarrow \ \Pr(X_1>0, X_2>0)=\frac{\alpha}{2\pi}.$$ Therefore if $M_x=\left[\begin{array}{cc}2x_1&-w\\-w&2x_2\end{array}\right]$ formula (\ref{FOR}) gives the following double integral on the domain $D=\{(t_1,t_2)
,\ \, t_1<x_1,\ t_2<x_2,\  w<2\sqrt{t_1t_2}\}$:
$$\arccos\frac{w}{2\sqrt{x_1x_2}}=\int_D
\frac{dt_1dt_2}{\sqrt{(x_1-t_1)(x_2-t_2)(4t_1t_2-w^2)}},$$ an identity not so easy to check directly.

\item Some comments about tentative  applications to Bayesian analysis  of $ MRIG_n$  are in order. Recall that a positive matrix $A=\rho I_n-C$ is called a $M$-matrix if $C=(c_{ij})_{1\leq i,j\leq n}$ is such that $c_{ij}\geq 0$ for all $i,j.$ Of course with our usual notation and for $x\in C_{W}$ then $M_x$ is a $M$-matrix: we have just to  define $c_{ij}=w_{ij}$ for $i\neq j$, $\rho=\max_i2x_i$ and $c_{ii}=\rho-2x_i$ for seeing this fact. The $M$-matrices are widely used in statistics since for $X\sim N(0,\Sigma)$ then the density $g(x)$ of $X$ has the $MTP_2$ property, namely for all $x,y\in \R^n$
$$g(\min(x_1,y_1)\ldots,\min(x_n,y_n))g(\max(x_1,y_1)\ldots,\max(x_n,y_n))\geq g(x)g(y)$$
if and only if $\Sigma^{-1}$ is a $M$-matrix: we refer for instance to Karlin and Rinott (1983) page 482 for this fact. In Theorem 3 of the same paper it is proved that for $X\sim N(0,\Sigma)$ and for all $i,j$ the covariance of $X_i,X_j$ conditioned by $\{X_k;\ 1\leq k\leq n,\  k\neq i,j\}$ is non negative if and only if $\Sigma^{-1}$ is a $M$-matrix.
From the point of view of Bayesian analysis two types of Gaussian models come to mind

\begin{enumerate}
\item $\{N(0, M^{-1}_{\theta})\ ; \ \theta\in C_W\}$. If $X\sim N(0, M^{-1}_{\theta})$
its density is $$\frac{1}{\sqrt{2\pi}^n}e^{-\frac{1}{2}x^*M_{\theta}x}\sqrt{\det M_{\theta}}$$
The densities have the $MTP_2$
property and the conditional covariances are all non negative. In order to use the $MRIG_n$ integral one  is tempted to consider the a priori measure
$$\pi(d\theta) =e^{-\frac{1}{2}b^*M_{\theta}^{-1}b}1_{C_W}(\theta)\frac{d\theta}{\det M_{\theta}}$$ which is unfortunately unbounded since $x\mapsto \int_{C_W}N(0,M_{\theta})(x)\pi(d\theta)$ is an unbounded density. From \eqref{STZI4} and the last comment before the proof of Theorem 2.2,  this density is proportional to $\prod_{i=1}^ne^{-a_i|x_i|}/|x_i|.$

\item  $\{N(0, M_{\theta})\ ; \ \theta\in C_W\}$. If $X\sim N(0, M_{\theta})$
its density is $$\frac{1}{\sqrt{2\pi}^n}e^{-\frac{1}{2}x^*M^{-1}_{\theta}x}\frac{1}{\sqrt{\det M_{\theta}}}$$ 
These densities have less attractive properties from the $MTP_2$ point of view. Nevertheless the a priori measure $$\pi(d\theta) =e^{-(a_1^2\theta_1+\cdots+a_n^2\theta_n)}1_{C_W}(\theta)d\theta$$ is bounded. However a major defect of this choice is the fact that $x\mapsto \int_{C_W}N(0,M^{-1}_{\theta})(x)\pi(d\theta)$  is computable (by \eqref{STZI4}) only if $x_1,\ldots,x_n$ are all non negative (again, see  example $n=2$ in Section 3)
\end{enumerate}\end{itemize}

\section{Continuity of the density of the $MRIG_n$ laws}
\noindent\textbf{Proposition 7.1.} If $b_1,\ldots,b_n\geq 0$  with $b\neq 0$ and if the graph associated to $W$ is connected then the function $$f(x)=e^{-\d b^*M_x^{-1}b}\frac{1}{\sqrt{\det M_x}}1_{C_W}(x)$$ is continuous on $\R^n.$

\vspace{4mm}\noindent\textbf{Proof.} The continuity of $f$ is clear outside of the boundary of $C_W$, namely outside of  the set $\partial C_W$ of $x\in \R^n$ such that $M_x$ is positive semidefinite with $\det M_x=0.$  In the sequel we fix $x\in \partial C_W$ and we prove the continuity of $f$ at this point $x.$

 \vspace{4mm}\noindent\textsc{First step.} We show that if $t=(t_1,\ldots,t_n)^*\in \R^n$ is  such that $M_xt=0$ and if $ t_{i_0}>0$ for some $i_0,$ then $t_i>0$ for all $i=1,\ldots,n.$ To see this, we use the notation $t_i^+=\max(0,t_i)$ , $ t_i^-=t_i^+-t_i$ and $$t^+=(t_1^+,\ldots,t_n^+)^*, \ \ t^-=t^+-t.$$ Since $0=M_xt=M_xt^+-M_xt^-$ we multiply by $(t^+)^*$ on the left for getting $(t^+)^*M_xt^+=(t^+)^*M_xt^-.$ Since $ t_i^+t_i^-=0$ we have that $(t^+)^*M_xt^-\leq 0.$ Since $M_x$ is positive semidefinite we have that $(t^+)^*M_xt^+= 0$ and therefore $M_xt^+= 0.$  Without loss of generality, assume that $t^+=(t_1,\ldots,t_k,0,\ldots,0)$ with $  t_1,\ldots,t_k>0$. Let us show that $k=n.$ Since $ t_{i_0}>0$ we have  $k>0.$ Suppose that $k<n$. We now split $M_x$ in blocks

$$M_x=\left[\begin{array}{cc}A&B\\B^*&C\end{array}\right]$$
where $A$ is a $(k,k)$ matrix. Clearly since $M_xt^+=0$ we get $B^*(t_1,\ldots,t_k)^*=0.$ Since it holds for all $t_j>0$ for $j=1,\ldots,k$ this implies that $B=0.$ This contradicts the fact that $G$ is connected, and finally $k=n.$

\vspace{4mm}\noindent\textsc{Second step.} We show that no principal minor of $M_x$ of order $n-1$ can be zero. Suppose for instance that the cofactor $C_{i_0}(x)$ of $2x_{i_0}$ is zero.
This implies that there exists a non-zero $t\in \R^n$ such that $M_xt=0$ and $ t_{i_0}=0.$ From the first step, this is impossible.

\vspace{4mm}\noindent\textsc{Third step.} Consider a sequence $(x_k)_{k=1}^{\infty}$ in $C_W$ converging to $x$ and let us show that $f(x_k)$ converges to zero. This is equivalent to show that

$$E_k= b^*M_{x_k}^{-1}b+\log \det M_{x_k}\to \infty$$
Recall that we have assumed that there exists $i_0$ such that  $b_{i_0}>0$. Recall also that all coefficients of $M_{x_k}^{-1}$ are non-negative.  As a consequence
$$E_k\geq b_{i_0}^2C_{i_0}(x_k)\frac{1}{\det M_{x_k}}+\log \det M_{x_k}$$ As  polynomials in $x_k$ we have that  $C_{i_0}(x_k)$ converges to $C_{i_0}(x)$ and $\det M_{x_k}$ converges to $\det M_x$. Since  $C_{i_0}(x_k)>0$ from the second step and since  $\det M_x=0,$ we have shown  that $E_k$ tends to infinity and the proof is done. $\square$

\section{Acknowlegments}We thank Christophe Sabot for having introduced the first author to   these questions and  for many useful comments. This research was partially supported for the second author  by the project 2016/21/B/ST1/00005 of the National Science Center, Poland. The first author is grateful for the hospitality of the   Faculty of Mathematics and Information Science of the Warsaw University of Technology.


\begin{thebibliography}{9}



\bibitem{BAPAT} Bapat, R.B.  (2010) \textit{Graphs  and Matrices}, Springer Universitext, Springer London Dordrecht Heidelberg NewYork.

 \bibitem{OLEBLAESILD}\textsc{Barndorff-Nielsen, O. E., Blaesild, P. and Seshadri, V. } (1992) Multivariate distributions with generalized inverse Gaussian marginals and associated Poisson mixtures. \textit{Can. J. Statist.} \textbf{20} , 109-120.


\bibitem{OLEKOUDOU} Barndorff-Nielsen, O. E. and Koudou, A. E.  (1998) Trees with random conductance and the (reciprocal) inverse Gaussian distribution. \textit{Adv. Appl. Probab.} \textbf{30} , 409-424.





\bibitem{OLERYDBERG}Barndorff-Nielsen, O. E. and Rydberg, T. H.  (2000) Exact distributional results for random resistance trees. \textit{Scand. J. Statist. } \textbf{27}(1) , 129-141.

\bibitem{BOBECKA} Bobecka, K. (2015) The Matsumoto-Yor property on trees for matrix variates of different dimensions. \textit{J. Multivar. Anal.} \textbf{141}, 22-44.


\bibitem{BOOLE}Boole, G. (1848) Th\'eor\`eme g\'en\'eral concernant l'int\'egration  d\'efinie. \textit{J.  Math. Pures  Appl. (1)},  \textbf{13}, 111-112.


\bibitem{DISERTORISZ10}Disertori, M.,  Spencer, T. and  Zinbauer, M. R.  (2010) Quasi Diffusion in a 3D Supersymmetry Hyperbolic Sigma Model, \textit{Comm. Math. Phys.} \textbf{300} 435-486.



\bibitem{DISERTORIMR17}Disertori, M.,  Merkl, F. and  Rolles, S. W. W. (2017) A supersymmetric approach to martingales related to the vertex-reinforced jump process, arXiv:1511.07157, \textit{ALEA Lat. J. Probab. Math. Stat.} 529-555.


\bibitem{KARLIN}Karlin, S. and Rinott, Y. (1983) $M$-matrices as covariance matrices of multinomial distributions. \textit{Linear Algebra and Its Applications} \textbf{52/53}, 419-438.





\bibitem{LETACW}Letac, G. and Weso\l owski, J. (2000) An independence property for the product of GIG and gamma laws. \textit{Ann. Probab.} \textbf{28}, 1371-1383.

\bibitem{MASSAMW04} Massam, H. and Weso\l owski, J. (2004) The Matsumoto-Yor property on trees. \textit{Bernoulli} \textbf{10(4)}, 685-700.

\bibitem{MASSAMW06} Massam, H. and Weso\l owski, J. (2006) The Matsumoto-Yor property and the structure of Wishart distributions. \textit{J. Multivariate Anal.} \textbf{97}, 103-123.

\bibitem{MYOR} Matsumoto, H. and Yor, M. (2001) An analogue of Pitman's $2M-X$ theorem for exponential Wiener functionals. Part II: The role of the generalized inverse Gaussian laws. \textit{Nagoya Math. J.} \textbf{162}, 65-68.

\bibitem{ST15} Sabot, C. and Tarr\`es, P.  (2015) Edge-reinforced random walk, vertex-reinforced jump process and the supersymmetric sigma model. \textit{J. Eur. Math. Soc.} \textbf{17} (9), 2353-2378.
\bibitem{SZ17}    Sabot, C.  and Zeng, X. (2017) Hitting times of interacting drifted Brownian motions and the vertex reinforced jump process, arXiv: 1704.05394.
\bibitem{SZ19}   Sabot, C.  and Zeng, X. (2019) A random Schr\"odinger operator associated with the Vertex Reinforced Jump Process on infinite graphs. \textit{J. Amer. Math. Soc.} \textbf{32}, 311-349. 
\bibitem {STZ17}Sabot, C., Tarr\`es, P.  and Zeng, X. (2017) The vertex reinforced jump process and a random Schrodinger operator on finite graphs.\textit{ Ann. Probab.} \textbf{45}, 3967-3986. 
\bibitem {SESHADRI}Seshadri, V.(1993) \textit{The Inverse Gaussian distribution} Clarendon Press, Oxford.


\bibitem{WATSON}Watson, G.N. (1966) \textit{A Treatise on the Theory of Bessel Functions} Cambridge University Press.

\end{thebibliography}
\end{document}